\documentclass[11pt,reqno,a4paper]{amsart}
\usepackage[margin=2.5cm,marginpar=2cm]{geometry}
\usepackage[all]{xy}
\usepackage{url}
\usepackage{amssymb}
\usepackage{hyperref}
\hypersetup{
  colorlinks,
  linkcolor={blue!50!black},
  citecolor={blue!50!black},
  urlcolor={blue!80!black}
}
\usepackage{bm}
\usepackage{tikz}
\usepackage{cite}
\usepackage{bbm}

\AtBeginDocument{%
\def\MR#1{}
}
\theoremstyle{plain}
\newtheorem{theorem}{Theorem}
\newtheorem{proposition}[theorem]{Proposition}

\newtheorem{conjecture}[theorem]{Conjecture}
\theoremstyle{definition}
\newtheorem{definition}[theorem]{Definition}
\newtheorem{example}[theorem]{Example}
\newtheorem{remark}[theorem]{Remark}
\newtheorem{question}[theorem]{Question}
\newtheorem{problem}[theorem]{Problem}
\newtheorem{convention}[theorem]{Convention}
\newtheorem{construction}[theorem]{Construction}
\newcommand{\PP}{\mathbb P}
\newcommand{\ZZ}{\mathbb Z}
\newcommand{\CC}{\mathbb C}
\newcommand{\QQ}{\mathbb Q}

\newcommand{\Aff}{\mathbb A}
\newcommand{\cG}{\mathcal{G}}
\newcommand{\cN}{\mathcal{N}}
\newcommand{\cM}{\mathcal{M}}
\newcommand{\cO}{\mathcal{O}}
\newcommand{\cF}{\mathcal{F}}
\newcommand{\cD}{\mathcal{D}}
\newcommand{\cK}{\mathcal{K}}
\newcommand{\tT}{\widetilde{T}}
\newcommand{\bx}{\bm{x}}
\newcommand{\abs}[1]{\left|#1\right|}
\newcommand{\hG}{\widehat{G}}
\newcommand{\one}{\mathbbm{1}}
\DeclareMathOperator{\coeff}{coeff}
\DeclareMathOperator{\Newt}{Newt}

\DeclareMathOperator{\Hom}{Hom}
\DeclareMathOperator{\GL}{GL}
\DeclareMathOperator{\SL}{SL}
\DeclareMathOperator{\PW}{PW}
\DeclareMathOperator{\Spec}{Spec}
\DeclareMathOperator{\Id}{Id}
\DeclareMathOperator{\Gr}{Gr}
\DeclareMathOperator{\Pic}{Pic}
\DeclareMathOperator{\conv}{conv}
\DeclareMathOperator{\Mon}{Mon}
\DeclareMathOperator{\verts}{vert}
\DeclareMathOperator{\crit}{crit}
\DeclareMathOperator{\im}{im}
\newcommand{\aff}{{\mathrm{aff}}}
\newcommand{\CY}{{Calabi\nobreakdash--Yau}}
\begin{document}
\author[A.\ M.\ Kasprzyk]{Alexander Kasprzyk}
\address{School of Mathematical Sciences\\University of Nottingham\\Nottingham\\UK}
\email{a.m.kasprzyk@nottingham.ac.uk}
\author[V.\ Przyjalkowski]{Victor Przyjalkowski}
\address{Steklov Mathematical Institute of Russian Academy of Sciences\\8 Gubkina street\\Moscow\\Russia}
\email{victorprz@mi-ras.ru, victorprz@gmail.com}
\title{Laurent polynomials in Mirror Symmetry: why and how?}
\thanks{AK was supported by EPSRC Fellowship~EP/N022513/1. VP was supported by President Grant MD-30.2020.1 and by the Foundation for the Advancement of Theoretical Physics and Mathematics ``BASIS''}
\begin{abstract}
We survey the approach to mirror symmetry via Laurent polynomials, outlining some of the main conjectures, problems, and questions related to the subject. We discuss: how to construct Landau--Ginzburg models for Fano varieties; how to apply them to classification problems; and how to compute invariants of Fano varieties via Landau--Ginzburg models.
\end{abstract}
\maketitle
\section{Introduction}
Mirror symmetry suggests a conjectural relationship between Fano varieties
and their Landau\nobreakdash--Ginzburg~(LG) models --- one dimensional families of \CY\ varieties dual to anticanonical sections of the Fano varieties. The duality relates symplectic properties of the Fano variety~$X$ (equipped with a complexification of a symplectic form on it) with algebraic properties of the dual LG~model~$(Y,w)$. The most general mirror symmetry conjectures (for example, those arising from Homological Mirror Symmetry) are hard to analyse. In this paper we discuss an effective approach to constructing LG~models of Fano varieties, and to computing their numerical properties, which confirms the general mirror symmetry expectation.

We are mostly interested in one of the two ``arrows'' of mirror symmetry that relate symplectic properties of a Fano variety~$X$ and algebraic properties of the dual LG~model~$w\colon Y\to \CC$. To consider~$X$ as a symplectic variety we first fix a symplectic form; in this paper we chose the anticanonical form (however many of the invariants we study do not depend on the choice of the form). We discuss two main problems: how to construct LG~models, and how to apply this to the problem of classifying Fano varieties; and how to compute invariants for an LG~model, and how to relate these invariants to the invariants of Fano varieties.
For the first problem we use an open chart (the algebraic torus) of the LG~model;
this enables us to apply the machinery of toric geometry and combinatorics. The second problem mostly deals with compactifications of these open charts and their cohomological invariants.

The exposition in this paper closely follows \cite{CCGGK13}~(for the material in~\S\ref{section:mirror partners}), 
\cite{Prz18b} (for the material in \S~\ref{section: toric degenerations} and \S~\ref{section: CY compactifications and toric LG}),
\cite{ACGK12}~(for the material in~\S\ref{sec:mutation}), \cite{CKPT21}~(for the material in~\S\ref{section: rigid MMLP}), \cite{HKP20,LP18}~(for the material in~\S\ref{section: KKP and P=W}), and \cite{ChP18}~(for the material in~\S\ref{section:anticanonical systems}).
\section{Two key examples}\label{section:two examples}
Many of the phenomena we wish to study already appear in the following two examples.

\begin{example}\label{example:P2}
Let~$X$ be the projective plane~$\PP^2$. To~$X$ we associate the formal power series
\[
\hG_X(t)=\sum_{k\geq 0}\frac{(3k)!}{(k!)^3}t^{3k}.
\]
Let us call the Laurent polynomial
\[
f=x+y+\frac{1}{xy}\in \CC[x^{\pm 1}, y^{\pm 1}]
\]
a mirror partner for~$X$. Denote the constant coefficient of~$f^k$ by~$\coeff_1(f^k)$ and set
\[
\pi_f(t)=\sum_{k\geq 0}\coeff_1(f^k)t^k.
\]
It is easy to see that~$\hG_X(t)=\pi_f(t)$.

Now consider the family~$\cF=\{ft=1\}$,~$t\in \PP^1\setminus\{0\}$, of fibres of the map~$(\CC^\times)^2\to \CC$ given by~$f$. This is a family of non-compact curves. Let us compactify it. For this consider the embeddings
\[
(\CC^\times)^2=\Spec\CC[x^{\pm 1}, y^{\pm 1}]\hookrightarrow \PP=\PP(x:y:z)
\]
and~$(\CC^\times)^2\times \left(\PP\setminus \{0\}\right)\hookrightarrow \PP\times \PP^1$, where the coordinates on~$\Aff^1$ and~$\PP^1$ are~$t$ and~$t_0$,~$t$, respectively. Let~$\widetilde Z$ be the closure of the graph of~$f$. The projection~$\widetilde{Z}\to \PP^1$ gives a structure of a rational elliptic surface. The variety~$\widetilde{Z}$ is singular; however it has du~Val singularities and admits a crepant resolution~$Z\to \widetilde{Z}$. The obtained elliptic surface has the fibre over~$\infty=(1:0)$ which is a wheel of nine smooth rational curves, and three other singular fibres having ordinary double points. Homological Mirror Symmetry conjecture for this family is studied in~\cite{AKO06}. Notice that~$\dim |{-K_X}|=9$.

Finally consider the Laurent polynomial
\[
f'=b+\frac{(a+1)^2}{ab^2}\in \CC[a^{\pm 1}, b^{\pm 1}].
\]
One can see that
\[
\pi_{f'}(t)=\sum_{k\ge 0}\frac{(3k)!}{(k!)^3}t^{3k}.
\]
This family of open curves which are fibres of~$f'$ can be compactified to the family of elliptic curves~$Z'\to \PP^1$. In fact one has~$Z'\cong Z$. The reason is that the families of fibres for~$f$ and~$f'$ are birational over the base; this can be seen using the birational transformation
\[
x=\frac{ab}{a+1},\qquad y=\frac{b}{a+1}.
\]
Note also that the Newton polytopes of~$f$ and~$f'$ --- that is, the convex hulls of exponents of monomials of~$f$ and~$f'$ --- are, respectively,
\[
\conv\{(1,0), (0,1), (-1,-1)\}\quad\text{ and }\quad
\conv\{(0,1),(1,-2),(-1,-2)\}.
\]
If we consider toric surfaces whose fans are generated by the cones over the faces of the Newton polytopes, we get, respectively,~$X=\PP^2$ and~$X'=\PP(1,1,4)$. Considering the second Veronese embedding of~$X'$ to~$\PP(1,1,1,2)$ with coordinates~$z_0$,~$z_1$,~$z_2$ of weight~$1$ and~$z_3$ of weight~$2$, one can describe~$X'$ as the quadric given by~$z_0z_1-z_2^2$. The projection of a general quadric in~$\PP(1,1,1,2)$ along the last coordinate gives the isomorphism of the general quadric and~$\PP^2$. Thus~$X$ degenerates to~$X'$.
\end{example}

\begin{example}
\label{example:C3}
Let~$X$ be a smooth cubic threefold. To~$X$ we associate the formal power series
\[
\hG_X(t)=\sum_{k\geq 0}\frac{(2k)!(3k)!}{(k!)^5}t^{2k}.
\]
Let us call the Laurent polynomial
\[
f=\frac{(x+y+1)^3}{xyz}+z\in \CC[x^{\pm 1}, y^{\pm 1},z^{\pm 1}]
\]
a mirror partner for~$X$. Denote again the constant coefficient of~$f^k$ by~$\coeff_1(f^k)$ and set
\[
\pi_f(t)=\sum_{k\geq 0}\coeff_1(f^k)t^k.
\]
It is easy to see that~$\hG_X(t)=\pi_f(t)$. Let~$\Delta\subset\cN\otimes \QQ$, where~$\cN=\ZZ^3$, be the convex hull of exponents of~$f$, and define
\[
\nabla=\{u\mid \langle u,v\rangle \geq -1\text{ for all } v\in\Delta\}\subset 
\Hom(\cN,\ZZ)\otimes\QQ
\]
to be the polytope dual to~$\Delta$. That is,~$\nabla$ is given by
\[
\conv\{(2,0,-1),(0,2,-1),(-2,-2,-1),(0,0,1)\}.
\]
Let~$T$ and~$T^\vee$ be the toric Fano varieties whose fans are generated by the cones spanned by the faces of, respectively,~$\Delta$ and~$\nabla$, so that~$T$ and~$T^\vee$ are dual toric varieties. Let~$\widetilde T^\vee$ be a toric variety whose rays are generated by the integral points on the boundary of~$\nabla$. One can check that~$\widetilde T^\vee$ is a crepant resolution of~$T^\vee$. Compactifying the family of fibres of~$f\colon(\CC^\times)^3\to \CC$ using the natural embedding~$(\CC^\times)^3\hookrightarrow \widetilde T^\vee$, we get a pencil of~K3 surfaces generated by its general element and the boundary divisor of~$\widetilde T^\vee$. One can check that after a resolution of the base locus of this family (by a sequence of blow ups of smooth curves) one arrives to the pencil of~K3 surfaces~$u\colon Z\to \PP^1$. One has~$-K_Z=u^{-1} (p)$, $p\in \PP^1$. The pencil has four singular fibres: two with ordinary double points; one (``over~$0$'') consisting of six smooth rational surfaces; and one (``over~$\infty$'') consisting of~$14$ smooth rational surfaces. The dual intersection complex of the latter fibre is homotopic to a two-dimensional sphere (so that it is a central fibre of Kulikov's type~III degeneration, see~\cite{Kul77}). Note that the choice of~$\widetilde T^\vee$ and the resolution of the base locus of the pencil is not unique, however all compactifications differ by flops, so the structure of the reducible fibre does not depend on the resolution. Note also that~$h^{12}(X)=5=6-1$ and that~$\dim |{-K_X}|=14$.
One can check that the Neron\nobreakdash--Severi lattice of general element of the family is
\[
M_3=H\oplus E_8(-1)\oplus E_{8}(-1)\oplus \langle -2\cdot 3\rangle,
\]
where~$H$ is a hyperbolic lattice, so that the pencil is a unique family of~K3 surfaces Dolgachev\nobreakdash--Nikulin dual to the family of anticanonical sections of~$X$ polarised by the generator of~$\Pic(X)$.

Since the degeneration at infinity is a Kulikov's type~III degeneration of~K3 surfaces, its monodromy is maximally unipotent. On the other hand, the monodromy at the fibre over~$0$ is quasiunipotent, but not unipotent. Note also that~$X$ is not rational.

Define
\[
f'=\frac{(a+b+1)^2}{abc}+c(a+b+1)\in \CC[a^{\pm 1}, b^{\pm 1},c^{\pm}].
\]
One can check that
\[
\pi_{f'}(t)=\sum_{k\geq 0}\frac{(2k)!(3k)!}{(k!)^5}t^{2k}.
\]
Proceeding as above, one can construct the pencil of~K3 surfaces~$u'\colon Z'\to \PP^1$. Moreover,~$Z$ and~$Z'$ differ by flops. This follows from the fact that general fibres of families given by~$f$ and~$f'$ are birational; the birational isomorphism is given by the change of variables
\[
x=a,\qquad y=b,\qquad z=c(a+b+1).
\]
Let~$\Delta'$ be the convex hull of exponents of~$f'$, and let~$T'$ be the toric Fano variety whose fan is generated by the cones spanned by the faces of~$\Delta'$. If the ambient four-dimensional projective space for~$X$ has coordinates~$z_0$,~$z_1$,~$z_2$,~$z_3$,~$z_4$, then~$T$ can be described as the toric cubic given by~$z_1z_2z_3-z_0^3$, whilst~$T'$ can be described as the toric cubic given by~$z_1z_2z_3-z_0^2z_4$. Thus,~$T$ and~$T'$ are degenerations of~$X$, and they deform one to each other.
\end{example}

In the following sections we generalise observations from Examples~\ref{example:P2} and~\ref{example:C3}, add additional observations, and discuss problems, questions, and conjectures related to the subject.

\section{Mirror partners}\label{section:mirror partners}
Let~$X$ be an~$n$-dimensional Fano variety. Let~$\one$ denote the fundamental class of~$X$ and let~$\cK\subset H_2(X,\ZZ)$ be the set of classes of effective curves.
The series
\[
\hG_X(t)=\sum_{a\in \ZZ_{\ge 0},\beta\in\cK} \left(-K_X\cdot \beta\right)!\langle\tau_{a}
{\one}\rangle_\beta \cdot t^{-K_X\cdot \beta}\in \CC[[t]],
\]
where~$\langle\tau_{a} \one\rangle_{\beta}$ is a~\emph{one-pointed genus~$0$ Gromov\nobreakdash--Witten invariant with descendants}, see~\cite[VI-2.1]{Ma99}, is called the~\emph{regularized quantum period} (or~\emph{a constant term of regularized~$I$-series}) for~$X$. We can write
\[
\hG_X(t)=1+\sum_{k=2}^\infty k!c_kt^k.
\]
Roughly speaking, the coefficients~$c_k$ encode the number of rational curves on~$X$ that pass through algebraic cycles on~$X$. The meaning of this series is the following. Homological Mirror Symmetry suggests that quantum cohomology is the Hochschild cohomology of the Fukaya category associated with~$X$ considered as a symplectic variety. Quantum multiplication (defined by three-pointed genus~$0$ prime Gromov\nobreakdash--Witten invariants) in this ring determines First and Second Dubrovin's connections. They give the quantum (regularized quantum, respectively)~$\cD$-module, and, according to~\cite{Gi96}, the series~$\hG_X$ is the solution of regularized quantum~$\cD$-module.

We come to our first question.

\begin{question}\label{question:reconstruct}
The regularized quantum period~$\hG_X$ is expected to characterise~$X$. Can this be proven? Is there an algorithmic way to reconstruct~$X$ from~$\hG_X$?
\end{question}

We expect that the regularized quantum~$\cD$-module is isomorphic to the Picard\nobreakdash--Fuchs one on the mirror side. The numerical evidence for this is the foundation of the mirror correspondence we consider. Moreover, we are looking for the dual to~$X$ as an algebraic torus~$(\CC^\times)^n$ with a complex-valued function. Choosing a basis we may assume that this function is represented by a Laurent polynomial, so we associate this polynomial to~$X$.

More precisely, let~$f\in\CC[x_1^{\pm1},\ldots,x_n^{\pm1}]$ be a Laurent polynomial. Associated to~$f$ is the~\emph{classical period}
\begin{equation}\label{eq:period}
\pi_f(t)=\left(\frac{1}{2\pi i}\right)^n\int_{\abs{x_1}=\ldots=\abs{x_n}=1}\frac{1}{1-tf}\frac{dx_1}{x_1}\cdots\frac{dx_n}{x_n},\qquad t\in\CC,\abs{t}\ll\infty.
\end{equation}
Expanding $\pi_f(t)$ as a series in $t$ one gets
\[
\pi_f(t)=\sum_{k=0}^\infty\coeff_1(f^k)t^k.
\]
Here~$\coeff_1(f^k)$ denotes the constant coefficient of~$f^k$. When
\[
\widehat{G}_X=\pi_f,
\]
mirror symmetry suggests that there is a close relationship between the geometry of~$X$ and~$f$, and we say that~$f$ is a~\emph{mirror partner} (or~\emph{weak LG~model}) for~$X$.

\begin{question}
Every Fano manifold in dimension~$n\leq 3$ has a mirror partner~\cite{Prz13, CCGK16}, many examples are known in dimension~$4$~\cite{CKP15,CGKS20}; in higher dimensions complete intersections in projective spaces~\cite{Prz13} and Grassmannians~\cite{PSh15b}
are proven to have mirror partners. Does this behaviour continue: that is, does every Fano manifold have a mirror partner?
\end{question}

Analogous to Question~\ref{question:reconstruct} above, we can ask:

\begin{question}\label{question:reconstruct_2}
Is there an algorithmic way to reconstruct~$X$ from a mirror partner~$f$? Note that a partial answer to Question~\ref{question:reconstruct_2} is given in~\cite{CKP19}, via the technique of ``Laurent inversion''.
\end{question}

\begin{example}[cf.\ Example~\ref{example:P2}]\label{eg:P2}
Consider~$\PP^2$. By Givental~\cite{Gi96} this has regularized quantum period
\[
\hG_{\PP^2}(t)=\sum_{k=0}^\infty\frac{(3k)!}{(k!)^3}t^{3k}=1+6t^3+90t^6+1680t^9+34650t^{12}+\cdots.
\]
The Laurent polynomial~$f=x+y+1/xy$ has classical period
\[
\pi_f(t)=\sum_{k=0}^\infty
{3k\choose k,k,k}
=\hG_{\PP^2}(t).
\]
Hence~$f$ is a mirror partner to~$\PP^2$.
\end{example}

The meaning of the classical period is that it is a period of the family of fibres of the map given by the Laurent polynomial; this period is given by taking residue of the form in the integral~\eqref{eq:period} and integrating it over the cycle whose~$S^1$-neighborhood is the standard~$n$\nobreakdash-cycle on the~$n$\nobreakdash-dimensional torus; see~\cite{Gol07,Prz08,CCGK16} for details. In other words, it is a solution for a Picard\nobreakdash--Fuchs differential operator for the family of fibres for the Laurent polynomial. We expect that this period is a period for family of fibrewise compactified fibres as well. The coefficients of the classical period are expected to be hypergeometric so, if one somehow knows ``enough'' terms of the expansion of the period, a recurrence relation on the period coefficients can be obtained (this is essentially linear algebra) and the Picard\nobreakdash--Fuchs operator derived. Alternatively, Lairez's generalised Griffiths--Dwork algorithm~\cite{La16} enables one to compute the Picard\nobreakdash--Fuchs operator directly from~$f$ with very high probability.

\begin{question}
What combinatorial or geometric properties of the Laurent polynomial give effective bounds on the number of terms of the classical period required to compute the Picard\nobreakdash--Fuchs differential operator?
\end{question}

\section{Toric degenerations}\label{section: toric degenerations}
Batyrev and Givental developed the mirror correspondence for toric varieties, and for complete intersections in a toric variety. Deformation invariance of quantum cohomology in smooth families suggests that the toric correspondence can be extended to the non-toric case via deformations to toric varieties. More precisely, we associate a (possibly singular) toric variety~$T_f$ to~$f$, given by taking the~\emph{spanning fan} (or~\emph{face fan}) of the Newton polytope~$\Delta=\Newt{f}\subset \cN\otimes_\ZZ\QQ$. That is, we take the fan in the lattice~$\cN\cong\ZZ^n$ whose cones span the faces of~$\Delta$. In the case of Example~\ref{eg:P2} we obtain the toric variety~$T_f=\PP^2$. In general we can not expect that if~$f$ is a mirror partner of a smooth Fano variety~$X$, then~$X$ degenerates to~$T_f$. Indeed, if~$f=f(x_1,\ldots,x_n)$ is a mirror partner to a Fano variety~$X$ and~$X$ degenerates to~$T_f$, then~$f'=f(x_1^2,x_2,\ldots,x_n)$ is a mirror partner of~$X$ again, while~$X$ does not degenerate to~$T_{f'}$. To avoid a subtlety about finite coverings of toric varieties, only Laurent polynomials~$f$ such that the exponents of monomials in~$f$ generate the lattice~$\cN$ should be considered.

\begin{example}
Consider the Laurent polynomials
\[
f=x+\frac{1}{x}+y+\frac{1}{y}\quad\text{ and }\quad g=xy+\frac{y}{x}+\frac{1}{xy}+\frac{x}{y}.
\]
Both polynomials have period sequence
\[
\pi_f(t)=\pi_g(t)=\sum_{k=0}^\infty\sum_{m=0}^k\frac{k!}{(m!)^2((k-m)!)^2}t^{2k}=1+4t^2+36t^4+400t^6+4900t^8+\cdots
\]
which coincides with the regularized quantum period of~$\PP^1\times\PP^1$. Indeed,~$T_f=\PP^1\times\PP^1$, however~$T_g=\PP^1\times\PP^1/(\ZZ/2)$ of anticanonical degree four. What has gone wrong here is that the exponents of monomials in~$g$ generate an index two sublattice in~$\cN$. In fact the ``correct'' Laurent polynomial supported on~$\Newt{g}$ is
\[
h=2x+xy+2y+\frac{y}{x}+\frac{2}{x}+\frac{1}{xy}+\frac{2}{y}+\frac{x}{y}.
\]
This has period sequence
\[
\pi_h(t)=1+20t^2+96t^3+1188t^4+10560t^5+111440t^6+\cdots
\]
and is seen to be a mirror partner for the smooth del~Pezzo surface~$X_{(2,2)}\subset\PP^4$ of anticanonical degree four.
\end{example}

We arrive at the following:

\begin{conjecture}\label{conjecture: toric degenerations}
Suppose that~$f$ is a mirror partner to~$X$. Then~$X$ admits a~\emph{$\QQ$\nobreakdash-Gorenstein (qG-) degeneration} to the singular toric variety~$T_f$.
\end{conjecture}

This conjecture has been studied and is supported in many cases: for example, del~Pezzo surfaces, Fano threefolds, and complete intersections~\cite{pragmatic,ILP13,KKPS19}. See also the beautiful three-dimensional example by Petracci~\cite{Pet20}.

In order for this construction to make sense, we assume that~$\Delta$ contains the origin in its strict interior. Note that this is not restrictive: if the origin is outside~$\Delta$ then the period of~$f$ must be a constant, and hence~$f$ cannot be a mirror partner to a Fano manifold; if the origin is contained in a proper face of~$\Delta$ then we can reduce to a lower-dimensional situation. We also require that the vertices of~$\Delta$ are primitive lattice points, thus ensuring that~$T_f$ is a toric Fano variety. That is,~$\Delta$ is a~\emph{Fano polytope} (see~\cite{KN13} for an overview of Fano polytopes).

\begin{remark}
Many of the combinatorial constructions described in this paper generalise if we relax the requirement that the vertices of~$\Delta$ are primitive. These polytopes correspond to toric Deligne\nobreakdash--Mumford stacks.
\end{remark}

It is important to emphasise that there is not a one-to-one correspondence between Laurent polynomial mirrors and Fano manifolds. Indeed, typically if there exists one mirror partner~$f$ for~$X$ then there exist infinitely many mirror partners\footnote{We mean this in a non-trivial way. Clearly monomial change of basis will always preserve the period sequence.}. Here the key process is~\emph{mutation}, which we describe in~\S\ref{sec:mutation}. This generates new Laurent polynomials with the same classical period.

\section{\CY\ compactifications and toric Landau\nobreakdash--Ginzburg models}\label{section: CY compactifications and toric LG}
Landau--Ginzburg~(LG) models are one-dimensional families of \CY\ varieties mirror dual to anticanonical sections of Fano varieties. Thus, these families should be proper. We expect that the proper families are compactifications of the mirror partners as families of hypersurfaces in tori. However, simply being a compactification of a mirror partner is not sufficient to be an appropriate mirror for Homological Mirror Symmetry. The important obstruction for this is that the compactified family should be a family of \CY\ varieties. This means that for the mirror partner~$f$ there should exist a commutative diagram
\[
\label{equation:CCGK-compactification}
\xymatrix{
(\mathbb C^\times)^n\ar@{^{(}->}[rr]\ar@{->}[d]_{f}&&Y\ar@{->}[d]^{w}\\
\CC\ar@{=}[rr]&&\CC}
\]
where~$w$ is proper, and~$Y$ is a smooth (open) \CY\ variety (so that fibres of~$w$ are also \CY). In this case we say that~$f$ satisfies the~\emph{\CY\ condition}. In the framework of recent interest to compactification of LG~models over infinity, we may strengthen this condition to the requirement of the existence of~\emph{log \CY\ compactification} via the extension of the diagram above to the commutative diagram
\[
\xymatrix{
(\mathbb C^\times)^n\ar@{^{(}->}[rr]\ar@{->}[d]_{f}&&Y\ar@{->}[d]^{w}\ar@{^{(}->}[rr]&&Z\ar@{->}[d]^{u}\\
\CC\ar@{=}[rr]&&\CC\ar@{^{(}->}[rr]&&\PP^1}
\]
where~$Z$ is a proper variety such that $u^{-1}(\infty)\sim -K_Z$.

There is no general procedure for the construction of a (log) \CY\ compactification; moreover, the existence may depend on particular coefficients of Laurent polynomial~$f$. However in many cases one can use the~\emph{log \CY\ compactification construction} described in~\cite{Prz17}.

\begin{construction}\label{construction:compactification}
Let~$\Delta$ be the Newton polytope of~$f$ and let~$\nabla$ be the dual polytope. Assume that~$\nabla$ is integral (that is,~$\Delta$ is~\emph{reflexive}), and let~$T^\vee$ be the toric variety given by the spanning fan of~$\nabla$. Assume that~$T^\vee$ admits a crepant resolution~$\tT^\vee\to T^\vee$. The family of fibres of~$f\colon (\CC^\times)^n\to \CC$ compactifies to an anticanonical pencil in~$T^\vee$. This compactified family $\cF$ is generated by its general member~$\cF_\lambda$ and the ``fibre over infinity''~$\cF_\infty$, which is nothing but the boundary divisor of~$\tT^\vee$. Finally, assume that the base set of the family~$\cF_\lambda\cap \cF_\infty$ is a union of smooth codimension two components (possibly with multiplicities). Blow these components up one-by-one to resolve the pencil. We obtain a family~$u\colon Z\to \PP^1$ such that~$Z$ is smooth and~$u^{-1}(\infty)\sim -K_Z$; this is the required log \CY\ compactification.
\end{construction}

Note that this procedure gives the description of the fibre~$u^{-1}(\infty)$, which is, up to codimension one, the boundary divisor of~$\tT^\vee$. It also gives the cohomology of~$Z$.

Construction~\ref{construction:compactification} is proven to be applicable (that is, all conditions are satisfied) for the Minkowski mirrors~\cite{ACGK12} for Fano threefolds, and for complete intersections; see~\cite{Prz17,Prz18a}. The main obstruction for the general case is that the polytope~$\Delta$ need not be reflexive. Fortunately, at least in some cases, the construction can be generalised.

\begin{example}[\!\!{\cite[Theorem 1.21]{Prz22}}]\label{example:singular LG}
Let~$X$ be a hypersurface of degree~$ad$ in
\[
\PP(\underbrace{1,\ldots,1}_{\alpha},d),
\]
where $\alpha=a(d-1)+1$,
with mirror partner
\[
f=\frac{(x_1+\cdots+x_{\alpha}+1)^{ad}}{x_1\cdots x_{\alpha}}.
\]
Let~$\Delta$ be the Newton polytope of~$f$. Compactify the family corresponding to~$f$ in the toric variety~$\tT^\vee$ defined by the spanning fan for the (non-integral) polytope~$\nabla=\Delta^\vee$; in fact~$\tT^\vee$ is a projective space. The support of the anticanonical divisor is the boundary divisor for~$\tT^\vee$; however the member of the family is~\emph{not} linearly equivalent to the anticanonical divisor, because the latter have multiplicity greater than one in one of the component of the boundary divisor of~$\tT^\vee$. After a carefully chosen resolution of the base locus for the compactified family, this component can be contracted to a (possibly singular) point; this gives a log \CY\ compactification. This compactification has a singular point over infinity, which cannot be avoided: there is no projective smooth log \CY\ compactification of the mirror partner.
\end{example}

\begin{problem}
Generalise Construction~\ref{construction:compactification} to the non-reflexive case.
\end{problem}

Of course, the output of Construction~\ref{construction:compactification} depends on a choice of a crepant resolution and a choice of the sequence of blow-ups. However, Hironaka-type arguments show that all (log) \CY\ compactifications of a given Laurent polynomial differ by flops. Moreover, (log) \CY\ compactifications of two mutation equivalent mirror partners also differ by flops.

\begin{question}\label{question:uniqueness of LG}
Is a (log) \CY\ compactification of a mirror partner to a given Fano variety uniquely defined up to flops?
\end{question}

A positive answer to Question~\ref{question:uniqueness of LG} is indicated by del~Pezzo surfaces and Picard rank one Fano threefolds. Indeed, if we assume that mirror partners for smooth del~Pezzo surfaces should be rigid maximally mutable Laurent polynomial (see~\S\ref{section: rigid MMLP}), then all such partners for a given surface are mutational equivalent, so their compactifications are isomorphic. The same argument works for Fano threefolds with very ample anticanonical class. Moreover, (log) \CY\ compactifications of Fano threefolds are families of~K3 surfaces. These~K3 surfaces are expected to be mirror dual to anticanonical sections of Fano threefolds. In particular, the~K3 surfaces are expected to be~\emph{Dolgachev\nobreakdash--Nikulin dual} (see~\cite{Do96}) to each other (cf.\ Example~\ref{example:C3}). Roughly speaking, Dolgachev\nobreakdash--Nikulin dual families of polarised~K3 surfaces are those whose algebraic and transcendental sublattices in the second cohomology lattice interchange. In particular, the general anticanonical section of a Picard rank one Fano threefold is polarised by a lattice of rank one. Thus its dual is polarised by a lattice of rank~$19$. Since there is a one-dimensional family of such~K3 surfaces, if we require the Dolgachev\nobreakdash--Nikulin duality for LG~models then the dual family is unique. The known mirrors of Picard rank one Fano threefolds satisfy Dolgachev\nobreakdash--Nikulin duality, see~\cite[5.4.3]{Prz18b}. Note that this argument does not assume that the LG~model we consider is a \CY\ compactification of a mirror partner.

The two requirements for mirror partners --- correspondence to toric degenerations and existence of \CY\ compactification --- give rise to the following definition.

\begin{definition}\label{definition: toric LG}
Let~$X$ be a smooth Fano variety of dimension~$n$. A Laurent polynomial~$f\in\CC[x_1^{\pm 1}, \ldots, x_n^{\pm 1}]$ is called a~\emph{toric Landau--Ginzburg model} of~$X$
if it satisfies the following three conditions.
\begin{description}
\item[Period condition] The polynomial~$f$ is a mirror partner for~$X$.
\item[\CY\ condition] The polynomial~$f$ satisfies the \CY\ condition.
\item[Toric condition] There is a flat degeneration~$X\rightsquigarrow T_f$.
\end{description}
\end{definition}

One can extend Definition~\ref{definition: toric LG} to the case of an arbitrary smooth projective variety~$X$. Moreover, if~$X$ is a smooth projective variety, and~$D$ is an element of the vector space~$\Pic(X)\otimes\CC$, one can define a toric LG~model for the pair~$(X,D)$ similarly to Definition~\ref{definition: toric LG}; see~\cite[Part~3]{Prz18b}.

\begin{question}
Construction~\ref{construction:compactification} suggests that a ``good'' mirror partner that satisfies
the toric condition (cf.~Conjecture~\ref{conjecture: toric degenerations}) will also satisfy the \CY\ condition. Is this always true? More generally, is it true that a mirror partner~$f$ for a Fano variety is a toric LG~model? What conditions on a Laurent polynomial guarantee this?
\end{question}

The following is the strong version of Mirror Symmetry of Variations of Hodge structures conjecture.

\begin{conjecture}[{see~\cite[Conjecture 38]{Prz13}}]\label{conjecture:MSVHS}
Any smooth Fano variety has a toric LG~model.
\end{conjecture}

Note that we associate a series~$\hG_X(t)$ and, thus, a toric LG~model, to a smooth Fano variety~$X$. However mirror symmetry associates an LG~model (as an algebraic variety) to~$X$ as a symplectic variety, or, in other words, to a pair of~$X$ and a divisor class on it. In fact we associate the series~$\hG_X(t)$ to a pair~$(X,-K_X)$. In the similar way, restricting the Gromov\nobreakdash--Witten series to an orbit of the torus~$\Hom(H^2(X,\ZZ),\CC^\times)$ to another orbit, generated by a class of another divisor~$D$, in an analogous way we can construct a series~$\hG_X^D(t)$, define a toric LG~model to it, and claim Conjecture~\ref{conjecture:MSVHS}; see~\cite[\S2]{PSh17}.

\section{Mutation}\label{sec:mutation}
Let~$f\in\CC[\bx^{\pm1},y^{\pm1}]$,~$F\in\CC[\bx^{\pm1}]$, where~$\bx=(x_1,\ldots,x_{n-1})$, and write
\[
f=\sum_{i\in\ZZ}P_i(\bx)y^i,\qquad\text{ where }P_i\in\CC[\bx^{\pm1}].
\]
Here all but finitely many of the~$P_i$ are zero. Suppose that there exists~$R_i\in\CC[\bx^{\pm1}]$ such that~$P_i=R_iF^{\abs{i}}$ for each~$i\in\ZZ_{<0}$. Define the map
\[
\begin{array}{r@{\ }c@{\ }l}
\mu\colon\CC(\bx,y)&\rightarrow&\CC(\bx,y)\\
\bx^ay^b&\mapsto&\bx^aF^by^b.
\end{array}
\]
Then we obtain a new Laurent polynomial
\[
g=\mu(f)=\sum_{i\in\ZZ_{<0}}R_iy^i+\sum_{j\in\ZZ_{\geq 0}}P_jF^jy^j\in\CC[\bx^{\pm1},y^{\pm1}].
\]
Furthermore, by an application of the change-of-variables formula to~\eqref{eq:period} we find that
\[
\pi_f(t)=\pi_g(t).
\]
Hence if~$f$ is a mirror partner to a Fano manifold~$X$, then~$g$ is also a mirror partner to~$X$.

\begin{definition}
Let~$\cN$ be a lattice of rank~$n$, and let~$w\in \cM=\Hom(\cN,\ZZ)$ be a primitive element in the dual lattice. Then~$w$ induces a grading on~$\CC[\cN]$. Let~$F\in\CC[w^\perp\cap \cN]$ be a Laurent polynomial in the zeroth graded piece of~$\CC[\cN]$, where
\[
w^\perp\cap \cN=\{v\in \cN\mid w(v)=0\}.
\]
The pair~$(w,F)$ defines an automorphism of~$\CC(\cN)$ via
\[
\begin{array}{r@{\ }c@{\ }l}
\mu_{(w,F)}\colon\CC(\cN)&\rightarrow&\CC(\cN)\\
\bx^v&\mapsto&\bx^vF^{w(v)}.
\end{array}
\]
We say that~$f\in\CC[\cN]$ is~\emph{mutable with respect to~$(w,F)$} if
\[
g=\mu_{(w,F)}(f)\in\CC[\cN].
\]
When this is the case, we call~$g$ a~\emph{mutation of~$f$}, and~$F$ a~\emph{mutation factor of~$f$}.
\end{definition}

For further details on mutation, see~\cite{ACGK12}. It was shown by Ilten~\cite{Il12} that if~$f$ and~$g$ are connected via a sequence of mutations, then~$T_f$ and~$T_g$ are related via qG-deformation; see also the generalisation~\cite{Pet21}.

\begin{example}[{see~\cite{GU10,AK16,KP12} and Example~\ref{example:P2}}]\label{eg:Markov}
Consider the Laurent polynomial~$f=x+y+1/xy$. Then $f$ is a mirror partner to~$\PP^2$. We can write
\[
f=\frac{1}{xy}(1+xy^2)+x.
\]
Taking~$w=(2,-1)\in \cM$ and~$F=1+xy^2$ we obtain the mutation
\[
g=\frac{1}{xy}+x(1+xy^2)^2
\]
where~$g$ is also a mirror partner to~$\PP^2$. As expected, the corresponding toric variety~$T_g=\PP(1,1,4)$ is a qG-deformation of~$\PP^2$, see Example~\ref{example:P2}.
We can continue mutating, obtaining a tree of qG-deformation equivalent toric varieties:
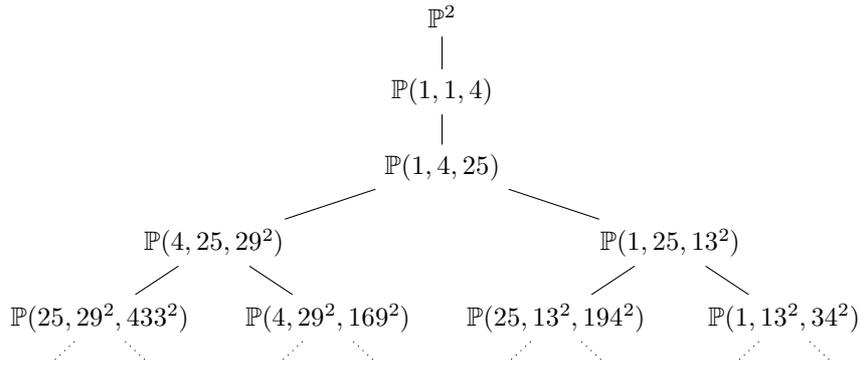
\begin{center}
\small
\begin{tikzpicture}[grow=down,level distance=1cm]
\tikzstyle{level 3}=[sibling distance=6cm]
\tikzstyle{level 4}=[sibling distance=3cm]
\tikzstyle{level 5}=[sibling distance=1.5cm,level distance=0.7cm]
\node {$\PP^2$}
child {node {$\PP(1,1,4)$}
child {node {$\PP(1,4,25)$}
child {node {$\PP(4,25,29^2)$}
child {node {$\PP(25,29^2,433^2)$}
child {node {} edge from parent[dotted]}
child {node {} edge from parent[dotted]}}
child {node {$\PP(4,29^2,169^2)$}
child {node {} edge from parent[dotted]}
child {node {} edge from parent[dotted]}}}
child {node {$\PP(1,25,13^2)$}
child {node {$\PP(25,13^2,194^2)$}
child {node {} edge from parent[dotted]}
child {node {} edge from parent[dotted]}}
child {node {$\PP(1,13^2,34^2)$}
child {node {} edge from parent[dotted]}
child {node {} edge from parent[dotted]}}}}
};
\end{tikzpicture}
\end{center}
Here the values~$(a,b,c)$ of the weighted projective space~$\PP(a^2,b^2,c^2)$ satisfy the~\emph{Markov equation}
\[
a^2+b^2+c^2=3abc,
\]
and are known as~\emph{Markov triples}. Mutation is equivalent, up to possible permutation of~$a$,~$b$, and~$c$, to the~\emph{Markov mutation}~$(a,b,c)\mapsto(3bc-a,b,c)$ of Markov triples. Note that by~\cite{HP10}, qG-degenerations of~$\PP^2$ are exactly~$\PP(a^2,b^2,c^2)$.
\end{example}

\begin{question}
The description of qG-deformations of~$\PP^2$ in Example~\ref{eg:Markov} is particularly elegant. Can a similar description be given for qG-deformations of~$\PP^1\times\PP^1$? Or, more generally, for each of the ten smooth del~Pezzo surfaces? Hacking\nobreakdash--Prokhorov~\cite{HP10} did this in the Picard rank one case, but what does this look like more generally?
\end{question}

\begin{conjecture}
Let~$f$ and~$g$ be mirror partners to~$X$. Then~$f$ and~$g$ are connected via a sequence of mutations.
\end{conjecture}

This, in particular, gives a uniqueness of rational LG~model of the same dimension as its dual smooth Fano variety.

\section{Rigid Maximally Mutable Laurent Polynomials}\label{section: rigid MMLP}
Any attempt at Fano classification via Laurent polynomials must address the following fundamental question.

\begin{question}
What class of Laurent polynomials are mirror partners to Fano varieties?
\end{question}

Here we have a conjectural answer: the~\emph{rigid maximally mutable Laurent polynomials (rigid MMLPs)} introduced in~\cite{CKPT21}. Roughly speaking, a Laurent polynomial is maximally mutable if it admits as many mutations as possible; an MMLP is rigid if it is uniquely determined by the mutations that it admits. We now make this definition precise. We begin by establishing some restrictions on the Laurent polynomials we consider.

\begin{convention}
A Laurent polynomial~$f\in\CC[\cN]$ is~\emph{normalised} if for all vertices~$v$ of~$\Newt{f}$, the coefficient of the monomial~$\bx^v$ in~$f$ is~$1$. We assume that all Laurent polynomials (and all mutation factors) from here onwards are normalised. Similarly, although our Laurent polynomials are defined over~$\CC$, our expectation is that after an appropriate choice a basis on the torus and scaling mirror partners have coefficients that are non-negative integers (although whether this assumption is correct remains an open question). We require that all Laurent polynomials (and all mutation factors) have non-negative integer coefficients. Furthermore, we require that for~$f \in \CC[\cN]$ the exponents of monomials in~$f$ generate~$\cN$.
\end{convention}

A transformation~$B\in\SL(\cN)$ is called a~\emph{$w$-shear}, where~$w\in \cM$, if~$B\mid_{w^\perp}=\Id$. Consider a mutation~$g=\mu_{(w,F)}(f)$. If we multiply the mutation factor~$F$ by a monomial~$\bx^v$,~$v\in w^\perp\cap \cN$, then~$\mu_{(w,F\bx^v)}(f)$ if related to~$g$ by a~$w$-shear. Thus considering~$F$ up to multiplication by monomials in~$\CC[w^\perp\cap \cN]$ gives~$g$ up to the action of~$w$-shears. Let
\[
\mu_{(w,F\bx^{w^\perp\cap \cN})}(f)
\]
denote the equivalence class given by~$w$-shears of~$g$.

Given a pair~$w\in \cM$,~$w$ primitive, and~$F\in\CC[w^\perp\cap \cN]$, write
\[
L(w,F)=\left(\langle w\rangle,F\bx^{w^\perp\cap \cN}\right).
\]
Here~$\langle w\rangle$ denotes the linear span of~$w$. We now define the mutation graph of a Laurent polynomial.

\begin{definition}
Let~$f\in\CC[\cN]$ be a Laurent polynomial. Define the graph~$G$ with vertices labelled by Laurent polynomials and edges labelled by pairs~$L(w,F)$ as follows. Write~$\ell(v)$ for the label of a vertex~$v$ of~$G$, and~$\ell(e)$ for the label of an edge~$e$ of~$G$.
\begin{enumerate}
\item
Begin with a vertex labelled by the Laurent polynomial~$f$.
\item
Given a vertex~$v$, set~$g=\ell(v)$. For each~$(w,F)$,~$\deg{F}>0$, such that~$g$ is mutable with respect to~$(w,F)$, and either:
\begin{enumerate}
\item
there does not exist an edge with endpoint~$v$ and label~$L(w,F)$; or
\item
for every edge~$e=\overline{vv'}$ with~$\ell(e)=L(w,F)$ we have that
\[
\ell(v')\not\in\mu_{(w,F\bx^{w^\perp\cap \cN})}(g)
\]
\end{enumerate}
pick a representative~$g'\in\mu_{(w,F\bx^{w^\perp\cap \cN})}(g)$ and add a new vertex~$v'$ and edge~$\overline{vv'}$ labelled by~$g'$ and~$L(w,F)$, respectively.
\end{enumerate}
The~\emph{mutation graph}~$\cG_f$ of~$f$ is defined by removing the labels from the edges of~$G$ and changing the labels of the vertices from~$G$ to the~$\GL(\cN)$-equivalence class of~$\Newt{g}$.
\end{definition}

We partially order the mutation graphs of Laurent polynomials: $\cG_f\prec\cG_g$ if there exists a label-preserving injection~$\cG_f\hookrightarrow\cG_g$.

\begin{definition}
A Laurent polynomial~$f$ is~\emph{maximally mutable} (or, for short,~$f$ is~\emph{MMLP}) if: $\Newt{f}$ is a Fano polytope; the constant term of~$f$ is zero; and~$\cG_f$ is maximal with respect to~$\prec$. A maximally mutable Laurent polynomial~$f$ is~\emph{rigid} if for all~$g$ such that the constant term of~$g$ is zero,~$\Newt{f}=\Newt{g}$, and~$\cG_f=\cG_g$, we have that~$f=g$.
\end{definition}

The close relationship between mutations of Laurent polynomials and cluster varieties suggests that being rigid should be a ``local'' property in the mutation graph. We have the following.

\begin{conjecture}
Let~$f\in\CC[\cN]$ be a Laurent polynomial such that~$\Newt{f}$ is a Fano polytope and the constant term of~$f$ is zero. Define
\[
S_f=\{(w,F)\mid f\text{ is mutable with respect to }(w,F)\}
\]
and, given any set~$S$ of pairs~$(w,F)$ with~$w\in \cM$,~$w$ primitive, and~$F\in\CC[w^\perp\cap \cN]$, define
\[
L_P(S)=\left\{f\in\CC[\cN]\, \left|\,
\parbox{9cm}{
$\Newt{f}=P$, the constant term of~$f$ is zero, and~$f$ is mutable with respect to~$(w,F)$ for all~$(w,F)\in S$
}
\right.\right\}.
\]
Then~$f$ is a rigid MMLP if and only if~$L_{\Newt{f}}(S_f)=\{f\}$.
\end{conjecture}

In two dimensions the picture is clear (see~\cite{pragmatic} for an overview). Fano polygons can be classified by their~\emph{singularity content}~\cite{AK14}. Those with singularity content~$(n,\varnothing)$ for some~$n$ fall into exactly ten mutation-equivalence classes~\cite[Theorem~1.2]{KNP17}, and each mutation class supports exactly one mutation class of rigid MMLPs~\cite[Theorem~3.9]{CKPT21}. These rigid MMLPs correspond one-to-one with qG-deformation families of smooth del~Pezzo surfaces. Under this correspondence, the classical period~$\pi_f$ of a rigid MMLP~$f$ matches with the regularized quantum period~$hG_X$ of the del~Pezzo surface~\cite[\S G]{CCGK16}. We obtain the following result.

\begin{theorem}[\!\!{\cite[Theorem~3.12]{CKPT21}}]
Mutation-equivalence classes of rigid MMLPs in two variables correspond one-to-one with qG-deformation families of smooth del~Pezzo surfaces.
\end{theorem}

A similar result holds in dimension three, building on the results of~\cite{ACGK12,CCGK16}:

\begin{theorem}[\!\!{\cite[Theorem~4.1]{CKPT21}}]
Mutation-equivalence classes of rigid MMLPs~$f$ such that $\Newt{f}$ is a three-dimensional reflexive polytopes correspond one-to-one to the~$98$ deformation families of three-dimensional Fano manifolds~$X$ with very ample~$-K_X$. Furthermore, each of the~$105$ deformation families of three-dimensional Fano manifolds has a rigid MMLP mirror.
\end{theorem}

\noindent
Furthermore, the four-dimensional mirrors in~\cite{CGKS20} have all been shown to be rigid MMLPs.

Notice that any simplicial terminal Fano polytope~$P\subset \cN\otimes_\ZZ\QQ$ supports a rigid MMLP
\[
f=\sum_{v\in\verts{P}}\bx^v.
\]
The variety $X_P$ is a~$\QQ$\nobreakdash-factorial Fano toric variety with at worst terminal singularities. Such varieties are know to be rigid under deformation~\cite{dFH12}. We are naturally led to consider a wider class of Fano varieties.

\begin{conjecture}[\!\!{\cite[Conjecture~5.1]{CKPT21}}]\label{conj:mmlps}
Rigid MMLPs in~$n$ variables (up to mutation) are in one-to-one correspondence with pairs~$(X,D)$, where~$X$ is a Fano~$n$-fold of class TG\footnote{A Fano variety~$X$ is of~\emph{class TG} if it admits a qG-degeneration with reduced fibres to a normal toric variety~\cite{pragmatic}.} with terminal locally toric qG-rigid singularities and~$D\in|{-K_X}|$ is a general element (up to qG-deformation). Under this correspondence, the classical period~$\pi_f$ of~$f$ agrees with the regularised quantum period~$\hG_X$ of~$X$, and~$X$ admits a qG-degeneration to the toric variety~$T_f$ given by the spanning fan of~$\Newt{f}$.
\end{conjecture}

\section{KKP and P=W conjectures}\label{section: KKP and P=W}
To approach the Katzarkov--Kontsevich--Pantev conjectures and Mirror P=W conjecture, we start from the following claim. Let~$X$ be a smooth Fano variety of dimension~$n$ and let~$Y$ be a \CY\ compactification of its toric LG~model. Set
\[
k_{Y}= \#\big(\text{irreducible components of all reducible fibres of~$Y$}\big)-\#\big(\text{reducible fibres}\big).
\]
Recall that the~\emph{primitive Hodge numbers} of~$X$ are defined as (see, for example,~\cite[p.~122]{GH78})
\[
h_{pr}^{p,q}(X)=
\left\{\begin{array}{ll}
h^{p,q}(X),&\text{ when }p\neq q;\\
h^{p,p}(X)-1,&\text{ when }p=q.
\end{array}
\right.
\]

\begin{conjecture}[\!\!{\cite[Conjecture 1.1]{PSh15}}]\label{conjecture: first Hodge}
$h_{pr}^{1,n-1}(X)=k_{Y}$.
\end{conjecture}

This conjecture is proven for certain toric LG~models of del~Pezzo surfaces~\cite[Proposition 1.20]{Prz22}, Fano threefolds~\cite[Main Theorem]{ChP18}, and complete intersections~\cite[Theorem 1.2]{PSh15}; see also~\cite{BGRS20}. Now we generalise Conjecture~\ref{conjecture: first Hodge} to other Hodge numbers.

It turns out that not only the reducible fibres of the LG~models themselves, but also the monodromy around them, affects the invariants of Fano varieties. The following result is given by comparing rationality of Picard rank one Fano threefolds studied by Iskovskikh and his school and Golyshev's computations of monodromies of their LG~models (cf.~Question~\ref{question:uniqueness of LG} and the discussion afterwards).

\begin{theorem}[\!\!{\cite[Theorem 3.3]{KP09}}]\label{theo:Gol-Isk}
Let~$X$ be a smooth Picard rank one Fano threefold whose compactified LG~model has a fibre with non-isolated singularities. Then the monodromy (in the second cohomology) at this fibre is unipotent if and only if~$X$ is rational.
\end{theorem}

\begin{problem}
Generalise Theorem~\ref{theo:Gol-Isk} to the higher Picard rank cases.
\end{problem}

The notion of log \CY\ compactification is very close to the notion of a tame compactified LG~model. We present it here in a reduced form adapted to our needs.

\begin{definition}[\!\!{\cite[Definition 2.4]{KKP17}}]\label{def-3}
A~\emph{tame compactified LG~model} is the data~$((Z,f),D_Z)$, where
\begin{enumerate}
\item~$Z$ is a smooth projective variety and~$f\colon Z\to \PP^1$ is a flat morphism.
\item~$D_Z=(\cup _i D^h_i)\cup (\cup _jD_j^v)$ is a reduced normal crossings divisor such that
\begin{itemize}
\item[(a)] $D^v=\cup _jD^v_j$ is a scheme-theoretical pole divisor of~$f$, i.e.~$f^{-1}(\infty)=D^v$. In particular~$ord _{D^v_j}(f)=-1$ for all~$j$;
\item[(b)] each component~$D_i^h$ of~$D^h=\cup _iD^h_i$ is smooth and horizontal for~$f$, i.e.~$f\vert _{D^h_i}$ is a flat morphism;
\item[(c)] the critical locus~$\crit(f)\subset Z$ does not intersect~$D^h$.
\end{itemize}
\item~$D_Z$ is an anticanonical divisor on~$Z$.
\end{enumerate}
One says that~$((Z,f),D_Z)$ is~\emph{a compactification of the LG~model}
$(Y,w)$ if in addition the following holds:
\begin{enumerate}\setcounter{enumi}{3}
\item~$Y=Z\setminus D_Z$,~$f\vert _Y=w$.
\end{enumerate}
\end{definition}

From now on we assume that $D^v=0$, so that $w$ is proper.
Note that the difference between log \CY\ compactifications of toric LG~models and tame compactified LG~models is that we allow the former to have singularities over infinity and do not require the fibre over infinity to be a normal crossing divisor. However the first issue does not effect statements about tame compactified LG~models, and the second does not appear in cases we know (say, in all cases when Construction~\ref{construction:compactification} is applicable).

In the following, all cohomology groups are taken with complex coefficients for the sake of simplicity. Let~$X$ be a Fano manifold. We assume that its mirror dual object~$(Y,w)$ admits a tame compactification. In~\cite{KKP17} Hodge-theoretic invariants~$f^{p,q}(Y,w)$ of an LG~model are constructed. We define
\[
f^{p,q}(Y,w) = \dim \Gr^F_qH^{p+q}(Y,V),
\]
where~$V$ is a smooth fibre of~$w$ and~$H^{p+q}(Y,V)$ is equipped with the natural mixed Hodge structure on the relative cohomology. This is equivalent to the definition in~\cite{KKP17,Ha17}.

\begin{conjecture}[\!\!{\cite[Conjecture 3.7]{KKP17}}]\label{conjecture: KKP mirror}
Let~$X$ and~$(Y,w)$ form a homological mirror pair. Then
\begin{equation}\label{eq:kkpeqq}
f^{p,q}(Y,w) = h^{\dim X - p,q}(X).
\end{equation}
\end{conjecture}

Conjecture~\ref{conjecture: KKP mirror} holds for \CY\ compactifications of certain toric LG~models of del~Pezzo surfaces~\cite[Theorem 12(ii)]{LP18} and Fano threefolds~\cite[Corollary]{ChP18}.

The motivation of the definition of another numbers that play a role of Hodge numbers for LG~models comes from Homological Mirror Symmetry, Hochschild homology identifications, and the identification of the monodromy operator with the Serre functor. Namely, assume that the LG~model~$(Y,w)$ is of Fano type (see~\cite[Definition 7]{LP18}) and is a mirror of a projective Fano manifold~$X$,~$\dim X=\dim Y$. Then by Homological Mirror Symmetry conjecture one expects an equivalence of categories
\begin{equation}\label{eq-of-cat}
D^b(coh\ X)\simeq FS((Y,w),\omega_Y),
\end{equation}
where~$D^b(coh\ X)$ is the bounded derived category of coherent sheaves on~$X$ and~$FS((Y,w),\omega_Y)$ is the Fukaya--Seidel category of the LG~model~$(Y,w)$ with an appropriate symplectic form~$\omega_Y$. This equivalence induces for each~$a$ an isomorphism of the Hochschild homology spaces
\[
HH_a(D^b(coh\ X))\simeq HH_a(FS((Y,w),\omega_Y)).
\]
It is known that
\begin{equation}\label{eq-of-hoch-derham}
HH_a(D^b(coh\ X))\simeq\bigoplus _{p-q=a}H^p(X,\Omega _X^q),
\end{equation}
and it is expected that
\begin{equation}\label{eq-of-hoch-fuk}
HH_a(FS((Y,w),\omega_Y))\simeq H^{n+a}(Y,V),
\end{equation}
where, as above,~$V$ is a smooth fibre of~$w$. The equivalence~\eqref{eq-of-cat} and isomorphisms~\eqref{eq-of-hoch-derham} and~\eqref{eq-of-hoch-fuk} suggest an isomorphism
\[
H^{n+a}(Y,V)=\bigoplus _{p-q=a}H^p(X,\Omega _X^q).
\]
Moreover, the equivalence~\eqref{eq-of-cat} identifies the Serre functors~$S_X$ and~$S_Y$ on the two categories. The functor~$S_X$ acts on the cohomology~$H^* (X)$, and the logarithm of this operator is equal (up to a sign) to the cup-product with~$c_1(K_X)$. Since~$X$ is Fano, the operator~$c_1(K_X)\cup (\,\cdot\,)$ is a Lefschetz operator on the space $\oplus_{p-q=a}H^p(X,\Omega _X^q)$ for each~$a$. On the other hand, the Serre functor~$S_Y$ induces an operator on the space~$H^{n+a}(Y,V)$ which is the inverse of the monodromy transformation~$M$ obtained by letting~$V$ vary in a small circle around~$\infty$. This suggests that the weight filtration for the nilpotent operator~$c_1(K_X)\cup (\,\cdot\,)$ on the space~$\bigoplus _{p-q=a}H^p(X,\Omega _X^q)$, that gives Hodge numbers for~$X$, should coincide with the similar filtration for the logarithm~$N$ of the operator~$M$ on~$H^{n+a}(Y,V)$. More precise,~$N$ gives a filtration~$\Mon$ on~$H^i(Y,V)$. We set
\[
h^{p,q}(Y,w) = \dim \Gr^{\Mon}_{p} H^{q}(Y,V).
\]

\begin{conjecture}[see {\cite[Conjecture 3.6]{KKP17}}]\label{conjecture:KKP h=f}
Let~$(Y,w)$ be a LG~model admitting a tame compactification. Then
\begin{equation}\label{eq:fanoht}
h^{p,q}(Y,w) = f^{p,q}(Y,w).
\end{equation}
\end{conjecture}

The following result is an extension of the main result of Shamoto~\cite{Sh17}.
\begin{proposition}[\!\!\cite{HKP20}]\label{corollary:lastones}
Let~$(Y,w)$ be a LG~model, and assume that~$w$ is a proper map. If~$H^i(Y)$ is Hodge\nobreakdash--Tate for all~$i$, then
\[
h^{p,q}(Y,w) = f^{p,q}(Y,w)
\]
for all~$p,q$.
\end{proposition}

Note that the output of the log \CY\ procedure is of Hodge\nobreakdash--Tate type provided the components of the blown up base locus are of Hodge\nobreakdash--Tate type as well, which usually holds. In particular, this holds for del~Pezzo surfaces, Fano threefolds and complete intersections.

Equalities~\eqref{eq:kkpeqq} and~\eqref{eq:fanoht} also hold for a smooth toric weak Fano threefolds~$X$ for which the map~$H^2(X) \rightarrow H^2(D)$ is injective for~$D$ a smooth anticanonical divisor, see~\cite{Ha17}.

Mirror symmetry constructions for Fano variety~$X$ consider it as both algebraic and symplectic variety. In other words, the input is an algebraic variety~$X$ equipped with a class of compexified symplectic form. If we do not mention it, then this symplectic form is anticanonical. However mirror duality can be strengthened: one can consider not a class of (anticanonical) divisors, but a certain simple normal crossing anticanonical divisor~$D$ on~$X$. We call such a pair~$(X,D)$~\emph{log \CY}. We abuse notation and call~$U = X \setminus D$ log \CY\ if the pair~$(X,D)$ is. One can classically equip an open \CY\ variety~$U=X\setminus D$ by a mixed Hodge structure. Ignoring the weight filtration on it, one can define Hodge numbers~$h^{p,q}(U)$. Mirror symmetry predicts that for the dual open \CY\ variety~$U^\vee$ one gets~$\mbox{$h^{p,q}(U)=h^{n-p,q}(U^\vee)$}$ for~$n=\dim U=\dim U^\vee$. The natural question is: can this prediction be extended to a duality of mixed Hodge structures, that is, to involve the weight filtration to the duality? The answer on this question is given by the~\emph{Mirror P=W conjecture}.

As we have mentioned above, if~$U$ and~$U^\vee$ are mirror log \CY\ manifolds, then we expect at first approximation that~$H^*_c(U)$ and~$H^*_c(U^\vee)$ are isomorphic as vector spaces (see~\cite[Table 1]{KKP17}) with different gradings. By Poincar\'e duality,~$H^i_c(U) \cong H^{\dim U - i}(U)$, hence we may equivalently deal with the cohomology rings of~$U$ and~$U^\vee$. Both~$H^*(U)$ and~$H^*(U^\vee)$ admit a mixed Hodge structure, which is composed of a decreasing Hodge filtration~$F^\bullet$ and an increasing weight filtration~$W_\bullet$. We define
\[
h^{p,q}(U) = \dim \Gr_F^q H^{p+q}(U).
\]
In analogy with classical mirror symmetry for compact \CY\ varieties, we might expect that if~$U$ and~$U^\vee$ are a homological mirror pair of log \CY\ manifolds of dimension~$n$, then
\begin{equation}\label{eq:logcyhn}
h^{p,q}(U) = h^{n - p,q}(U^\vee).
\end{equation}
This seems to be true --- it is checked in many cases in~\cite{HKP20} --- but it ignores the weight filtration in cohomology. It would be desirable to determine whether the weight filtration on~$H^*(U)$ is reflected by a filtration on the cohomology of~$U^\vee$. The first step in this is to remark that the geometry of~$D = X \setminus U$ and the residues of holomorphic forms on~$X$ with log poles along~$D$ can be used to determine the weight filtration~$W_\bullet$ on~$H^*(U)$. The weight filtration depends on the existence of a projective simple normal crossings compactification~$X$ of~$U$, but is independent of the choice of compactification, hence it is a canonical invariant of~$U$. So if a mirror dual filtration exists, is plausible that it can be constructed via information dual to that of the components of~$D$ but be independent the choice of~$D$.

Starting with a log \CY\ manifold~$U$ and a simple normal crossings compactification~$X$ of~$U$ with~$D = X \setminus U$, each irreducible component~$D_i$,~$i=1,\ldots,k$, of~$D$ determines a regular function~$w_i$ on the mirror~$U^\vee$, see~\cite{aur1,aur2,aak}. Therefore, if there is a filtration on~$H^*(U^\vee)$ dual to the weight filtration on~$H^*(U)$, it should be determined by the functions~$w_{1},\ldots , w_{k}$.

There are several possible filtrations on cohomology that can be constructed from~$(w_1,\dots, w_k)$, but the most relevant seems to be the~\emph{flag filtration} \cite{dcm}, which is defined as follows. Let~$w$ denote the map~$(w_1,\dots, w_k)\colon U^\vee\rightarrow \CC^k$. Choose a generic flag of linear sub\-spa\-ces
$$\Lambda_k \subset \Lambda_{k-1} \subset \dots \subset \Lambda_0 =\CC^k$$
so that~$\dim \Lambda_i = k-i$ and let~$U^\vee_i = w^{-1}(\Lambda_i)$. Then, for any coefficient ring~$R$, the flag filtration on~$H^*(U^\vee;R)$ is defined as\footnote{Note that this agrees with the definition of~\cite{dcm} up to a shift by~$j$.}
\[
P_{r}H^j(U^\vee;R) = \ker(H^j(U^\vee;R) \longrightarrow H^j(U^\vee_{r+1};R)).
\]
According to de Cataldo and Migliorini~\cite{dcm}, if~$w$ is proper, then~$P_\bullet$ can be identified with the~\emph{perverse Leray filtration} of the map~$w$, hence it only depends on the map~$w$. In all of the cases that we know of, the maps~$w_{1},\dots, w_{k}$ generate~$\CC[U^\vee]$, in which case the map~$w\colon U^\vee \rightarrow\im(U^\vee)$ is the affinization map of~$U^\vee$, hence, in these cases at least,~$P_\bullet$ is~\emph{intrinsic to~$U^\vee$} and does not depend on our original choice of~$w_{1},\dots, w_{k}$. Thus we have two filtrations, which are built from data which correspond to one another under mirror symmetry, and which are intrinsic to~$U$ and~$U^\vee$ respectively.

\begin{definition}[\!\!{\cite[Definition 1.1]{HKP20}}]
Consider a quasiprojective variety~$M$ over the complex numbers~$\CC$ and assume that the affinization map~$f^\aff\colon M \rightarrow\Spec(\CC[M])$ is proper. We define the~\emph{perverse mixed Hodge polynomial} of a quasiprojective variety~$M$ to be
\[
\PW_{M}(u,t,w,p) = \sum_{a,b,r,s} (\dim \Gr_F^a \Gr^W_{s+b}\Gr^P_{r} (H^{s}(M)))u^a t^{s}w^{b}p^r,
\]
where~$P_\bullet$ is the flag filtration taken with respect to~$f^\aff$ and~$W$ denotes the~$\CC$-linear extension of the weight filtration.
\end{definition}

The following is called~\emph{Mirror P=W conjecture}.

\begin{conjecture}[\!\!{\cite[Conjecture 1.2]{HKP20}}]\label{con:filt}
Let~$U$ be a log \CY\ variety and assume that its homological mirror~$U^\vee$ is also a log \CY\ variety whose dimension is the same as that of~$U$. Let~$n = \dim U = \dim U^\vee$. Then
\[
\PW_{U}(u^{-1}t^{-2}, t,p, w)u^nt^n = \PW_{U^\vee}(u,t,w,p).
\]
\end{conjecture}

\begin{theorem}[\!\!\cite{HKP20}]
Let~$(X,D)$ be a pair consisting of smooth Fano surface or threefold~$X$ and a smooth anticanonical divisor~$D$ on it. Let~$(Y,w)$ be its compactified LG~model constructed in~\cite{AKO06} and~\cite{Prz17}. Then Conjecture~\ref{con:filt} holds for them.
\end{theorem}

When the map~$w$ is proper we expect that~$X$ admits a smooth anticanonical divisor~$D$ so that~$X \setminus D$ and~$Y$ form a homological mirror pair. Therefore, equality~\eqref{eq:logcyhn} should hold between~$Y$ and~$X \setminus D$. Furthermore, we expect that a general smooth fibre~$V$ of~$w$ is \CY\, and is the homological mirror of~$D$, so we expect that
\begin{equation}\label{eq:cptcyhn}
h^{p,q}(V) = h^{\dim{X} -1 -p,q}(D).
\end{equation}
Conjecture~\ref{con:filt} links equality~\eqref{eq:kkpeqq} with equalities~\eqref{eq:logcyhn} and~\eqref{eq:cptcyhn}.

\begin{theorem}[\!\!\cite{HKP20}]
Let~$X$ be a projective manifold with a smooth anticanonical divisor~$D$ in it, and let~$U = X\setminus D$. Let~$(Y,w)$ be a LG~model so that~$w$ is proper and let~$V$ be a smooth fibre of~$w$. If~$V$ and~$D$ satisfy~\eqref{eq:cptcyhn}, and~$Y$ and~$U$ satisfy equality~\eqref{eq:logcyhn}, then Conjecture~\ref{con:filt} implies equality~\eqref{eq:kkpeqq}.
\end{theorem}

\begin{theorem}[\cite{HKP20}]
Let~$X$ be a Fano manifold with a smooth anticanonical divisor~$D$ in~$X$, and~$(Y,w)$ be a LG~model so that~$w$ is proper. Assume that Conjecture~\ref{con:filt} holds between~$Y$ and~$U = X \setminus D$. Then
\[
f^{p,q}(Y,w) = h^{p,q}(Y,w).
\]
\end{theorem}

\section{Anticanonical linear systems}\label{section:anticanonical systems}
Conjecture~\ref{conjecture: first Hodge} relates the number of components of LG~models with the Hodge number of the corresponding Fano varieties. 
Katzarkov--Kontsevich--Pantev conjectures generalise this conjecture to other Hodge numbers. More precisely,~\cite[Theorem 4.8]{Ha17} says that for LG~model~$(Y,w)$ for a Fano threefold the equality~$f^{1,1}(Y,w)=k_Y$ holds if the LG~model satisfies certain natural conditions provided by Construction~\ref{construction:compactification}. Thus, Construction~\ref{construction:compactification} plays an important role in numerical conjectures for Fano--LG correspondence. Let us discuss some other implications.

Let~$X$ be a smooth Fano threefold and let~$f$ be its toric LG~model. Let~$\Delta$ be a Newton polytope for~$f$. Consider the flat degeneration of~$X$ to the toric Fano variety~$T_f$ whose spanning polytope is~$\Delta$. Since this degeneration is flat, one has
\[
\chi\left(\cO_X(-K_X)\right)=\chi\left(\cO_{T_f}(-K_X)\right).
\]
On the other hand, since~$T_f$ is toric, its singularities are Kawamata log terminal by~\cite[Proposition~3.7]{Ko95}. Applying Kodaira vanishing (see, for example,~\cite[Theorem~2.70]{KM98}) to~$X$ and~$T_f$, one gets
\[
h^i(\cO_X(-K_{X}))=h^i(\cO_{T_f}(-K_{T_f}))=0
\]
for~$i>0$, so that
\[
h^0(-K_X)=h^0(-K_{T_f}).
\]
The anticanonical linear system of~$T_f$ can be described as a linear system of Laurent polynomials supported on the dual polytope~$\nabla$, see, for instance,~\cite[\S 6.3]{Da78}. Suppose that Construction~\ref{construction:compactification} is applicable for~$f$. In particular,~$\nabla$ is integral and~$T^\vee$ admits a toric crepant resolution~$\widetilde T^\vee\to T^\vee$. The dimension of the anticanonical linear system of~$T_f$ is the number of integral points on the boundary of~$\nabla$. Since these boundary points are in one-to-one correspondence with boundary divisors of~$\widetilde T^\vee$ and, thus, with irreducible components of the fibre~$u^{-1}(\infty)$. This motivates the following conjecture.

\begin{conjecture}[\!\!{\cite[Conjecture 1.6]{ChP20}}]\label{conjecture:linear systems}
Let~$X$ be a smooth Fano variety, and let~$(Z,u)$ be a log \CY\ compactification of its toric LG~model~$f$. Then the fibre~$u^{-1}(\infty)$ consists of
\[
\chi\big(\cO_X(-K_X)\big)-1=h^0\big(\cO_X(-K_{X})\big)-1
\]
irreducible components.
\end{conjecture}

As we have mentioned, this conjecture is proved for rigid maximally-mutable toric LG~models of smooth Fano threefolds and for Givental's toric LG~models of ``good'' toric Fano varieties, see~\cite[Theorem~1.7]{ChP20}. Moreover, from Example~\ref{example:singular LG} one can see that Conjecture~\ref{conjecture:linear systems} holds even in the cases when Construction~\ref{construction:compactification} is not applicable.

\begin{remark}\label{remark:non-vanishing}
Let us notice that Conjecture~\ref{conjecture:linear systems} together with Conjecture~\ref{conjecture:MSVHS} imply that
\[
h^0\big(\cO_X(-K_{X})\big)\ge 2,
\]
which is only known for~$\dim(X)\le 5$ (see \cite[Theorem~1.7]{HV11},~\cite[Theorem~1.1.1]{HS19}). Note also that Kawamata's \cite[Conjecture~2.1]{Ka00} implies that~$h^0(\cO_X(-K_{X}))\ge 1$.
\end{remark}

Homological Mirror Symmetry conjecture suggests that the monodromy around~$u^{-1}(\infty)$ is maximally unipotent (see~\cite[\S2.2]{KKP17}). Thus, if the fibre~$u^{-1}(\infty)$ is a divisor with simple normal crossing singularities, then its dual intersection complex is expected to be homeomorphic to a sphere of dimension~$n-1$ (see~\cite[Question~7]{KoXu16}). This follows from~\cite[Proposition~8]{KoXu16} for~$\mbox{$n\le 5$}$. However, we cannot always expect~$u^{-1}(\infty)$ to be a divisor with simple normal crossing singularities. The example is a toric LG~model~$(Y,w)$ a smooth intersection of two general sextics in~$\PP(1,1,1,2,2,3,3)$, see~\cite[Example 1.9]{ChP20}. On the other hand, if we take a log resolution of the pair~$(Z,u^{-1}(\infty))$, that is, if we blow up~$Z$ to make the fibre over infinity a normal crossing divisor, then the dual intersection graph is homeomorphic to a~$3$-dimensional sphere, so we can expect that the answer on analogue of~\cite[Question~7]{KoXu16} holds.

\begin{problem}
Define a dual intersection complex for the degenerations of \CY\ varieties.
\end{problem}

We expect that, at least for LG~models, the answer on~\cite[Question~7]{KoXu16} for this definition is positive, cf.~Example~\ref{example:singular LG}.



\begin{thebibliography}{10}

\bibitem{aak}
Mohammed Abouzaid, Denis Auroux, and Ludmil Katzarkov.
\newblock Lagrangian fibrations on blowups of toric varieties and mirror
  symmetry for hypersurfaces.
\newblock {\em Publ. Math. Inst. Hautes \'{E}tudes Sci.}, 123:199--282, 2016.

\bibitem{pragmatic}
Mohammad Akhtar, Tom Coates, Alessio Corti, Liana Heuberger, Alexander~M.
  Kasprzyk, Alessandro Oneto, Andrea Petracci, Thomas Prince, and Ketil
  Tveiten.
\newblock Mirror symmetry and the classification of orbifold del {P}ezzo
  surfaces.
\newblock {\em Proc. Amer. Math. Soc.}, 144(2):513--527, 2016.

\bibitem{ACGK12}
Mohammad Akhtar, Tom Coates, Sergey Galkin, and Alexander~M. Kasprzyk.
\newblock Minkowski polynomials and mutations.
\newblock {\em SIGMA Symmetry Integrability Geom. Methods Appl.}, 8:Paper 094,
  17, 2012.

\bibitem{AK14}
Mohammad Akhtar and Alexander~M. Kasprzyk.
\newblock Singularity content.
\newblock \href{http://arxiv.org/abs/1401.5458}{\texttt{arXiv:1401.5458
  [math.AG]}}, 2014.

\bibitem{AK16}
Mohammad Akhtar and Alexander~M. Kasprzyk.
\newblock Mutations of fake weighted projective planes.
\newblock {\em Proc. Edinb. Math. Soc. (2)}, 59(2):271--285, 2016.

\bibitem{aur1}
Denis Auroux.
\newblock Mirror symmetry and {$T$}-duality in the complement of an
  anticanonical divisor.
\newblock {\em J. G\"{o}kova Geom. Topol. GGT}, 1:51--91, 2007.

\bibitem{aur2}
Denis Auroux.
\newblock Special {L}agrangian fibrations, wall-crossing, and mirror symmetry.
\newblock In {\em Surveys in differential geometry. {V}ol. {XIII}. {G}eometry,
  analysis, and algebraic geometry: forty years of the {J}ournal of
  {D}ifferential {G}eometry}, volume~13 of {\em Surv. Differ. Geom.}, pages
  1--47. Int. Press, Somerville, MA, 2009.

\bibitem{AKO06}
Denis Auroux, Ludmil Katzarkov, and Dmitri Orlov.
\newblock Mirror symmetry for del {P}ezzo surfaces: vanishing cycles and
  coherent sheaves.
\newblock {\em Invent. Math.}, 166(3):537--582, 2006.

\bibitem{BGRS20}
Edoardo Ballico, Elizabeth Gasparim, Francisco Rubilar, and Luiz~A.B.
  San~Martin.
\newblock {KKP} conjecture for minimal adjoint orbits.
\newblock \href{https://arxiv.org/abs/1901.07939}{\texttt{arXiv:1901.07939
  [math.AG]}}, 2019.

\bibitem{ChP18}
Ivan Cheltsov and Victor Przyjalkowski.
\newblock {K}atzarkov--{K}ontsevich--{P}antev conjecture for {F}ano threefolds.
\newblock \href{https://arxiv.org/abs/1809.09218}{\texttt{arXiv:1809.09218
  [math.AG]}}, 2018.

\bibitem{ChP20}
Ivan Cheltsov and Victor Przyjalkowski.
\newblock Fibers over infinity of {L}andau--{G}inzburg models.
\newblock \href{https://arxiv.org/abs/2005.01534}{\texttt{arXiv:2005.01534
  [math.AG]}}, 2020.

\bibitem{CCGGK13}
Tom Coates, Alessio Corti, Sergey Galkin, Vasily Golyshev, and Alexander~M.
  Kasprzyk.
\newblock Mirror symmetry and {F}ano manifolds.
\newblock In {\em European {C}ongress of {M}athematics}, pages 285--300. Eur.
  Math. Soc., Z\"{u}rich, 2013.

\bibitem{CCGK16}
Tom Coates, Alessio Corti, Sergey Galkin, and Alexander~M. Kasprzyk.
\newblock Quantum periods for 3-dimensional {F}ano manifolds.
\newblock {\em Geom. Topol.}, 20(1):103--256, 2016.

\bibitem{CGKS20}
Tom Coates, Sergey Galkin, Alexander~M. Kasprzyk, and Andrew Strangeway.
\newblock Quantum periods for certain four-dimensional {F}ano manifolds.
\newblock {\em Exp. Math.}, 29(2):183--221, 2020.

\bibitem{CKPT21}
Tom Coates, Alexander~M. Kasprzyk, Giuseppe Pitton, and Ketil Tveiten.
\newblock Maximally mutable {L}aurent polynomials.
\newblock {\em Proc. A.}, 477(2254):Paper No. 20210584, 21, 2021.

\bibitem{CKP15}
Tom Coates, Alexander~M. Kasprzyk, and Thomas Prince.
\newblock Four-dimensional {F}ano toric complete intersections.
\newblock {\em Proc. A.}, 471(2175):20140704, 14, 2015.

\bibitem{CKP19}
Tom Coates, Alexander~M. Kasprzyk, and Thomas Prince.
\newblock Laurent inversion.
\newblock {\em Pure Appl. Math. Q.}, 15(4):1135--1179, 2019.

\bibitem{Da78}
Vladimir~I. Danilov.
\newblock The geometry of toric varieties.
\newblock {\em Uspekhi Mat. Nauk}, 33(2(200)):85--134, 247, 1978.

\bibitem{dcm}
Mark Andrea~A. de~Cataldo and Luca Migliorini.
\newblock The perverse filtration and the {L}efschetz hyperplane theorem.
\newblock {\em Ann. of Math. (2)}, 171(3):2089--2113, 2010.

\bibitem{dFH12}
Tommaso de~Fernex and Christopher~D. Hacon.
\newblock Rigidity properties of {F}ano varieties.
\newblock In {\em Current developments in algebraic geometry}, volume~59 of
  {\em Math. Sci. Res. Inst. Publ.}, pages 113--127. Cambridge Univ. Press,
  Cambridge, 2012.

\bibitem{Do96}
Igor~V. Dolgachev.
\newblock Mirror symmetry for lattice polarized {$K3$} surfaces.
\newblock volume~81, pages 2599--2630. 1996.
\newblock Algebraic geometry, 4.

\bibitem{GU10}
Sergey Galkin and Alexandr Usnich.
\newblock Laurent phenomenon for {G}inzburg--{L}andau potential.
\newblock \href{http://db.ipmu.jp/ipmu/sysimg/ipmu/417.pdf}{IPMU preprint
  10-0100}, 2010.

\bibitem{Gi96}
Alexander~B. Givental.
\newblock Equivariant {G}romov--{W}itten invariants.
\newblock {\em Internat. Math. Res. Notices}, (13):613--663, 1996.

\bibitem{Gol07}
Vasily~V. Golyshev.
\newblock Classification problems and mirror duality.
\newblock In {\em Surveys in geometry and number theory: reports on
  contemporary {R}ussian mathematics}, volume 338 of {\em London Math. Soc.
  Lecture Note Ser.}, pages 88--121. Cambridge Univ. Press, Cambridge, 2007.

\bibitem{GH78}
Phillip Griffiths and Joseph Harris.
\newblock {\em Principles of algebraic geometry}.
\newblock Pure and Applied Mathematics. Wiley-Interscience [John Wiley \&
  Sons], New York, 1978.

\bibitem{HP10}
Paul Hacking and Yuri Prokhorov.
\newblock Smoothable del {P}ezzo surfaces with quotient singularities.
\newblock {\em Compos. Math.}, 146(1):169--192, 2010.

\bibitem{Ha17}
Andrew Harder.
\newblock Hodge numbers of {L}andau--{G}inzburg models.
\newblock {\em Adv. Math.}, 378:Paper No. 107436, 40, 2021.

\bibitem{HS19}
Andreas H\"{o}ring and Robert \'{S}miech.
\newblock Anticanonical system of {F}ano fivefolds.
\newblock {\em Math. Nachr.}, 293(1):115--119, 2020.

\bibitem{HV11}
Andreas H\"{o}ring and Claire Voisin.
\newblock Anticanonical divisors and curve classes on {F}ano manifolds.
\newblock {\em Pure Appl. Math. Q.}, 7(4, Special Issue: In memory of Eckart
  Viehweg):1371--1393, 2011.

\bibitem{Il12}
Nathan~Owen Ilten.
\newblock Mutations of {L}aurent polynomials and flat families with toric
  fibers.
\newblock {\em SIGMA Symmetry Integrability Geom. Methods Appl.}, 8:Paper 047,
  7, 2012.

\bibitem{ILP13}
Nathan~Owen Ilten, Jacob Lewis, and Victor Przyjalkowski.
\newblock Toric degenerations of {F}ano threefolds giving weak
  {L}andau--{G}inzburg models.
\newblock {\em J. Algebra}, 374:104--121, 2013.

\bibitem{KKPS19}
Alexander~M. Kasprzyk, Ludmil Katzarkov, Victor Przyjalkowski, and Dmitrijs
  Sakovics.
\newblock Projecting {F}anos in the mirror.
\newblock \href{https://arxiv.org/abs/1904.02194}{\texttt{arXiv:1904.02194
  [math.AG]}}, 2019.

\bibitem{KN13}
Alexander~M. Kasprzyk and Benjamin Nill.
\newblock Fano polytopes.
\newblock In {\em Strings, gauge fields, and the geometry behind}, pages
  349--364. World Sci. Publ., Hackensack, NJ, 2013.

\bibitem{KNP17}
Alexander~M. Kasprzyk, Benjamin Nill, and Thomas Prince.
\newblock Minimality and mutation-equivalence of polygons.
\newblock {\em Forum Math. Sigma}, 5:Paper No. e18, 48, 2017.

\bibitem{KKP17}
Ludmil Katzarkov, Maxim Kontsevich, and Tony Pantev.
\newblock Bogomolov--{T}ian--{T}odorov theorems for {L}andau--{G}inzburg
  models.
\newblock {\em J. Differential Geom.}, 105(1):55--117, 2017.

\bibitem{KP09}
Ludmil Katzarkov and Victor Przyjalkowski.
\newblock Generalized homological mirror symmetry and cubics.
\newblock {\em Tr. Mat. Inst. Steklova}, 264(Mnogomernaya Algebraicheskaya
  Geometriya):94--102, 2009.

\bibitem{KP12}
Ludmil Katzarkov and Victor Przyjalkowski.
\newblock Landau--{G}inzburg models---old and new.
\newblock In {\em Proceedings of the {G}\"{o}kova {G}eometry-{T}opology
  {C}onference 2011}, pages 97--124. Int. Press, Somerville, MA, 2012.

\bibitem{HKP20}
Ludmil Katzarkov, Victor Przyjalkowski, and Andrew Harder.
\newblock {$\textrm{P}=\textrm{W}$} {P}henomena.
\newblock {\em Mat. Zametki}, 108(1):33--46, 2020.

\bibitem{Ka00}
Yujiro Kawamata.
\newblock On effective non-vanishing and base-point-freeness.
\newblock volume~4, pages 173--181. 2000.
\newblock Kodaira's issue.

\bibitem{Ko95}
J\'{a}nos Koll\'{a}r.
\newblock Singularities of pairs.
\newblock In {\em Algebraic geometry---{S}anta {C}ruz 1995}, volume~62 of {\em
  Proc. Sympos. Pure Math.}, pages 221--287. Amer. Math. Soc., Providence, RI,
  1997.

\bibitem{KM98}
J\'{a}nos Koll\'{a}r and Shigefumi Mori.
\newblock {\em Birational geometry of algebraic varieties}, volume 134 of {\em
  Cambridge Tracts in Mathematics}.
\newblock Cambridge University Press, Cambridge, 1998.
\newblock With the collaboration of C. H. Clemens and A. Corti, Translated from
  the 1998 Japanese original.

\bibitem{KoXu16}
J\'{a}nos Koll\'{a}r and Chenyang Xu.
\newblock The dual complex of {C}alabi--{Y}au pairs.
\newblock {\em Invent. Math.}, 205(3):527--557, 2016.

\bibitem{Kul77}
Viktor~S. Kulikov.
\newblock Degenerations of {$K3$} surfaces and {E}nriques surfaces.
\newblock {\em Izv. Akad. Nauk SSSR Ser. Mat.}, 41(5):1008--1042, 1199, 1977.

\bibitem{La16}
Pierre Lairez.
\newblock Computing periods of rational integrals.
\newblock {\em Math. Comp.}, 85(300):1719--1752, 2016.

\bibitem{LP18}
Valery Lunts and Victor Przyjalkowski.
\newblock Landau--{G}inzburg {H}odge numbers for mirrors of del {P}ezzo
  surfaces.
\newblock {\em Adv. Math.}, 329:189--216, 2018.

\bibitem{Ma99}
Yuri~I. Manin.
\newblock {\em Frobenius manifolds, quantum cohomology, and moduli spaces},
  volume~47 of {\em American Mathematical Society Colloquium Publications}.
\newblock American Mathematical Society, Providence, RI, 1999.

\bibitem{Pet20}
Andrea Petracci.
\newblock An example of mirror symmetry for {F}ano threefolds.
\newblock In {\em Birational geometry and moduli spaces}, volume~39 of {\em
  Springer INdAM Ser.}, pages 173--188. Springer, Cham, 2020.

\bibitem{Pet21}
Andrea Petracci.
\newblock Homogeneous deformations of toric pairs.
\newblock {\em Manuscripta Math.}, 166(1-2):37--72, 2021.

\bibitem{Prz08}
Victor Przyjalkowski.
\newblock On {L}andau--{G}inzburg models for {F}ano varieties.
\newblock {\em Commun. Number Theory Phys.}, 1(4):713--728, 2007.

\bibitem{Prz13}
Victor Przyjalkowski.
\newblock Weak {L}andau--{G}inzburg models of smooth {F}ano threefolds.
\newblock {\em Izv. Ross. Akad. Nauk Ser. Mat.}, 77(4):135--160, 2013.

\bibitem{Prz17}
Victor Przyjalkowski.
\newblock Calabi--{Y}au compactifications of toric {L}andau--{G}inzburg models
  for smooth {F}ano threefolds.
\newblock {\em Mat. Sb.}, 208(7):84--108, 2017.

\bibitem{Prz18a}
Victor Przyjalkowski.
\newblock On the {C}alabi--{Y}au compactifications of toric
  {L}andau--{G}inzburg models for {F}ano complete intersections.
\newblock {\em Mat. Zametki}, 103(1):111--119, 2018.

\bibitem{Prz18b}
Victor Przyjalkowski.
\newblock Toric {L}andau--{G}inzburg models.
\newblock {\em Uspekhi Mat. Nauk}, 73(6(444)):95--190, 2018.

\bibitem{Prz22}
Victor Przyjalkowski.
\newblock On singular log {C}alabi--{Y}au compactifications of
  {L}andau--{G}inzburg models.
\newblock \href{https://arxiv.org/abs/2102.01388}{\texttt{arXiv:2102.01388
  [math.AG]}}, 2021.

\bibitem{PSh15b}
Victor Przyjalkowski and Constantin Shramov.
\newblock Laurent phenomenon for {L}andau-{G}inzburg models of complete
  intersections in {G}rassmannians.
\newblock {\em Proc. Steklov Inst. Math.}, 290(1):91--102, 2015.
\newblock Published in Russian in Tr. Mat. Inst. Steklova {{\bf{2}}90} (2015),
  102--113.

\bibitem{PSh15}
Victor Przyjalkowski and Constantin Shramov.
\newblock On {H}odge numbers of complete intersections and {L}andau--{G}inzburg
  models.
\newblock {\em Int. Math. Res. Not. IMRN}, (21):11302--11332, 2015.

\bibitem{PSh17}
Victor Przyjalkowski and Constantin Shramov.
\newblock Laurent phenomenon for {L}andau--{G}inzburg models of complete
  intersections in {G}rassmannians of planes.
\newblock {\em Bull. Korean Math. Soc.}, 54(5):1527--1575, 2017.

\bibitem{Sh17}
Yota Shamoto.
\newblock Hodge--{T}ate conditions for {L}andau--{G}inzburg models.
\newblock {\em Publ. Res. Inst. Math. Sci.}, 54(3):469--515, 2018.

\end{thebibliography}
\end{document}